\documentclass[12pt]{article}

\usepackage{amssymb,amsmath}

\newtheorem{thm}{Theorem}[section]
\newtheorem{lem}{Lemma}[section]
\newtheorem{C}{Corollary}[section]

\newtheorem{dfn}{Definition}[section]

\numberwithin{equation}{section}

\def\qed{\hfill $\square$}
\def\comp{\leavevmode\raise.2ex\hbox{${\scriptstyle\mathchar"020E}$}}
\def\circ{\comp}
\def\P{\mathbf{P}}
\def\E{\mathbf{E}}
\def\supp{{{\rm supp}\,}}
\def\proof{\noindent{\it Proof.\ \ }}
\def\SS{\mathbf{S}}

\def\BB{\mathcal{B}}

\def\DD{\mathbf{D}}
\def\FF{\mathcal{F}}
\def\EE{\mathcal{E}}
\def\R{\mathbb{R}}
\def\ss{{{\mathtt s}}}

\def\ang #1{{\langle #1\rangle}}

\begin{document}

\title {\bf The uniqueness of symmetrizing measure and linear diffusions}
\author{Xing Fang, Jiangang Ying
\footnote{The research of this author is supported in part by NSFC Grant No. 10671036}, Minzhi Zhao
\footnote{The research of this author is supported in part by NSFC Grant No. 10601047}}

\date{}
 \maketitle
\bigskip

\noindent{\bf 2000 MR subject classification} 60J45, 60J65

\noindent {\bf Key words.} symmetrizing measure, linear diffusion, Dirichlet space, regular
subspace

\begin{abstract}
In this short article, we shall study one-dimensional local Dirichlet spaces. One result, which has
its independent interest, is to prove that irreducibility implies the uniqueness of symmetrizing
measure for right Markov processes. The other result is to give a representation for any 1-dim
local, irreducible and regular Dirichlet space and a necessary and sufficient condition for a
Dirichlet space to be regular subspace of another Dirichlet space.
\end{abstract}

\section{Introduction}

Due to the pioneering works of Feller, one-dimensional diffusion has been a mature and very
interesting topic in theory of Markov processes with its simplicity and clarity. There are a lot of
literatures on this topic, e.g., Ito-Mckean\cite{im}, Revuz-Yor\cite{ry}, Rogers-Williams
\cite{rw}, among those most influential. As we shall see, one-dimensional irreducible diffusion is
always symmetric. Thus it has no loss of generality that Dirichlet form approach is introduced to
investigate one-dimensional diffusions. In this article, we shall discuss the properties of
Dirichlet spaces associated with one-dimensional diffusions, and study one-dimensional diffusions
by means of Dirichlet forms. At first a representation of the Dirichlet form associated with a
one-dimensional diffusion will be formulated since we have not seen it explicitly in literature.
%This representation makes it legitimate to study

\section{The uniqueness of symmetrizing measure}

We first present a theorem which states a condition for uniqueness of symmetrizing measure and will
be used later. This kind of results may be known in some other forms. We begin with a general right
Markov process $X=(X_t,\P^x)$ on state space $E$ with semigroup $(P_t)$ and resolvent
$(U^{\alpha})$. It is easy to see from the right continuity that for $x\in E$ and a finely open
subset $D$, $\P^x(T_D<\infty)>0$ if and only if $U^{\alpha}1_D(x)>0$. The process $X$ is called
irreducible if $\P^x(T_D<\infty)>0$ for any $x\in E$ and a finely open subset $D$, where $T_D$ is
the hitting time of $D$.

\begin{lem} \label{l:2.1} The following statements are equivalent.
\begin{enumerate}
\item $X$ is irreducible.
\item $U^\alpha 1_D$ is positive everywhere on $E$ for any non-empty finely open set $D$.
\item $U^\alpha 1_A$ is either identically zero or positive everywhere on $E$ for any Borel set $A$
or, in other words, $\{U^{\alpha}(x,\cdot):x\in E\}$ are all mutually absolutely continuous. .
\item All non-trivial excessive measures are mutually absolutely
continuous.
\end{enumerate}
\end{lem}

\proof The equivalence of (1) and (2) is easy. We shall prove that they are equivalent to (3). We
may assume $\alpha =0$. Suppose (1) is true.  If $U 1_A$ is not identically zero, then there exists
$\delta>0$ such that $D:=\{U 1_A>\delta\}$ is non-empty. Since $U 1_A$ is excessive and thus finely
continuous, $D$ is finely open and the fine closure of $D$ is contained in $\{U 1_A\ge \delta\}$.
Then $$U 1_A(x)\ge P_D U 1_A(x)=\E^x\left( U 1_A(X_{T_D})\right)\ge \delta \P^x(T_D<\infty)>0.$$
Conversely suppose (3) is true. Then for any finely open set $D$, by the right continuity of $X$,
$U 1_D(x)>0$ for any $x\in D$. Therefore $U 1_D$ is positive everywhere on $E$.

Let $\xi$ be an excessive measure. Since $\alpha \xi U^\alpha\le \xi$, $\xi(A)=0$ implies that $\xi
U^{\alpha}(A)=0$. However $\xi$ is non-trivial. Thus it follows from (3) that $U^{\alpha}1_A\equiv
0$, i.e., $A$ is potential zero. Conversely if $A$ is potential zero, then $\xi(A)=0$ for any
excessive measure $\xi$. Therefore (3) implies (4).

Assume (4) holds. Since $U^{\alpha}(x,\cdot)$ is excessive for all $x$ and hence they are
equivalent. This implies (3). \qed

A Borel set $A$ is called of potential zero if $U^{\alpha} 1_A$ is identically zero for some
$\alpha\ge 0$ (thus for all $\alpha\ge 0$). A $\sigma$-finite measure $\mu$ on $E$ is said to be a
symmetrizing measure of $X$ or $X$ is said to be $\mu$-symmetric if
$$(P_t u,v)_{\mu}=(u,P_t v)_{\mu}$$ for any measurable $u,v\ge 0$ and $t>0$. It is easy to check
that any symmetrizing measure is excessive and an excessive measure does not charge any set of
potential zero.

\begin{thm} \label{t:2.1}
Assume that $X$ is irreducible. Then the symmetrizing measure of $X$ is unique up to a constant.
More precisely if both $\mu$ and $\nu$ are non-trivial symmetrizing measures of $X$, then
$\nu=c\mu$ with a positive constant $c$.
\end{thm}

\proof First of all there exists a measurable set $H$ such that both $\mu(H)$ and $\nu(H)$ are
positive and finite, because $\mu$ and $\nu$ are equivalent by Lemma~\ref{l:2.1}. This is actually
true when both measures are $\sigma$-finite and one is absolutely continuous with respect to
another. Indeed, assume that $\nu\ll\mu$. Since $\nu$ is non-trivial and $\sigma$-finite, we may
find a measurable set $B$ such that $0<\nu(B)<\infty$. Then $\mu(B)>0$. % and $\nu(H)<\infty$.
Since $\mu$ is $\sigma$-finite, there exist $A_n\uparrow E$ such that $0<\mu(A_n)<\infty$. Then $\nu(A_n\cap B)\uparrow \nu(B)$ and $\mu(A_n\cap B)\uparrow \mu(B)$. Hence there exists some $n$ such that $\nu(A_n\cap B)>0$. Take $H=A_n\cap B$, which makes both $\mu(H)$ and $\nu(H)$ positive and finite. %If $\nu(H)=0$, we may replace $\nu$
%by $\nu+\mu$. Then the following argument will show that $\mu$ and $\nu+\mu$ are linear, and so are
%$\mu$ and $\nu$.

Set $c=\nu(H)/\mu(H)$. We may assume that $c=1$ without loss of generality.
 Let $m=\mu+\nu$.
Then there is $f_1,f_2\ge 0$ such $\mu=f_1\cdot m$ and $\nu=f_2\cdot m$.  Let $A=\{f_1>f_2\}$,
$B=\{f_1=f_2\}$ and $C=\{f_1<f_2\}$.

We shall show that $\nu=\mu$. Otherwise $\mu(A)>0$ or $\nu(C)>0$. We assume that $\mu(A)>0$ without
loss of generality. Since $\mu$ is $\sigma$-finite, there is $A_n\in\BB(E)$ such that $A_n\subseteq
A$, $\mu(A_n)<\infty$ and $A_n\uparrow A$. Let $D=B\cup C$. For any integer $n$ and $\alpha>0$,
$$ (U^\alpha 1_{A_n}, 1_D)_\mu\le (U^\alpha 1_{A_n},1_D)_\nu=(U^\alpha
1_D,1_{A_n})_\nu\le (U^\alpha 1_D, 1_{A_n})_\mu.$$ Since $(U^\alpha 1_{A_n}, 1_D)_\mu= (U^\alpha
1_D, 1_{A_n})_\mu$, it follows that $(U^\alpha 1_D,1_{A_n})_\nu= (U^\alpha 1_D, 1_{A_n})_\mu$. Thus
we have
$$(U^\alpha 1_D,(1-\frac {f_2}{f_1})1_{A_n})_\mu=(U^\alpha 1_D,1_{A_n})_\mu-
(U^\alpha 1_D, 1_{A_n})_\nu=0.$$ Since $1-\frac{f_2}{f_1}>0$ on $A$, let $n$ go to infinity and by
the monotone convergence theorem we get that $(U^\alpha 1_D, 1_A)_\mu=0$. The irreducibility of $X$
implies that $U^\alpha 1_D=0$ identically or $D$ is of potential zero. Therefore
$$\mu(D)=\nu(D)=0.$$
%The equality $$(U^\alpha 1_D, 1)_\mu=(U^\alpha 1,
%1_D)_\mu=0$$ leads to that $U^\alpha 1_D=0$ $\mu$-a.e. Since $U^\alpha 1_D$ is
%lower-semi-continuous and $\supp\mu=E$, $U^\alpha(x,D)=0$ for all $x\in E$. Therefore $\nu(D)=0$ by
%the equality $$(U^\alpha 1,1_D)_\nu=(U^\alpha 1_D,1)_\nu=0,$$ and $\mu(D)\le\nu(D)=0$.
Consequently, $$0=\mu(H)-\nu(H)=\int_{H\cap A}(1-\frac {f_2}{f_1})d\mu$$ which leads to that
$\mu(H\cap A)=0$ and also $\mu(H)=0$. The contradiction implies that $\nu=\mu$.
  \qed
%\section{Representation of one-dimensional local Dirichlet space}
\medskip

The following example shows that the condition that any point may reach any finely open set is
needed. Actually we may easily see that it is also necessary in the sense that if $X$ has a unique
symmetrizing measure $m$, then $X$, restricted on the fine support of $m$, is irreducible.

 \noindent \bf Example: \rm Let $J={1\over
4}(\delta_1+\delta_{-1}+\delta_{\sqrt2}+\delta_{-\sqrt2})$ defined on $\mathbb{R}$ and
$\pi=\{\pi_t\}_{t>0}$ the corresponding symmetric convolution semigroup; i.e.,
$\hat\pi_t(x)=e^{-t\phi(x)}$ with $$\phi(x)=\int (1-\cos xy)J(dy) ={1\over 2}(1-\cos x)+{1\over
2}(1-\cos {\sqrt2x}).$$ Let $N=\{n+m{\sqrt2}\ :\ n,m\ {\rm are\ integers}\}$ and $\mu=\sum_{x\in
N}\delta_x.$ Then $\mu$ is $\sigma$-finite and also a symmetrizing measure. It is easy to check
that any point may reach any open set but not any finely open set.
\\

%\end{document}

It is known that the fine topology is determined by the process and hard to identify usually. Hence
it is hard to verify sometimes the irreducibility defined in the theorem. However under LSC,
namely, assuming that $U^{\alpha}1_B$ is lower-semi-continuous for any Borel subset $B$ of $E$, the
irreducibility is equivalent to the weaker one, which is easier to verify: $\P^x(T_D<\infty)>0$ for
any $x\in E$ and open subset $D\subset E$.
\medskip

\noindent \bf Remark \rm As a remark, we would like to present a slight more general result which
was provided by Masatoshi Fukushima in his comment to this theorem.

Suppose that $X$ is $\mu$-symmetric. The following two definitions refer to Definition 2.1.1
\cite{CF}. A Borel subset $A$ is called $(P_t)$-invariant if $1_A\cdot P_t (1_{A^c}f)=0$ a.e. $\mu$
for all $t>0$ and $f\in L^2(E,\mu)$, and $X$ is $\mu$-irreducible if any $(P_t)$-invariant set is
trivial in the sense that either $\mu(A)=0$ or $\mu(A^c)=0$. Then the following statements are
equivalent due to Theorem 3.5.6\cite{CF} and a similar proof of Lemma \ref{l:2.1}.
\begin{itemize}
\item[(1)] $X$ is $\mu$-irreducible;
\item[(2)] If $D$ is finely open and $\mu(D)>0$, then $\P^x(T_D<\infty)>0$ for q.e. $x\in E$;
\item[(3)] $U^\alpha1_D>0$ q.e. for any finely open $D$ with $\mu(D)>0$;
\item[(4)] $U^\alpha 1_A$ is either 0 q.e. or positive q.e. for every Borel subset $A$.
\end{itemize}
It follows that if $X$ is $\mu$-irreducible, then all non-trivial excessive measures charging no
$\mu$-polar sets are equivalent. Hence following the proof of Theorem~\ref{t:2.1}, we have its
Fukushima's version.

\begin{thm} Assume that a Borel right process $X$ is $\mu$-irreducible with respect to some
non-trivial symmetrizing measure $\mu$ of $X$. If $\nu$ is a symmetrizing measure of $X$ charging
no $\mu$-polar sets, then $\nu=c\cdot \mu$ for some constant $c\ge 0$.
\end{thm}

\section{Dirichlet forms on intervals}

Let $I$ be an interval or a connected subset of $\R$ and $I^{\circ}$ its interior. Denote by
$\SS(I)$ the totality of strictly increasing continuous functions on $I$. Let $\ss\in \SS(I)$. Let
$m$ and $k$ two Radon measures on $I$ with $\supp(m)=I$. Define a symmetric form
$(\EE^{(\ss,m,k)},\FF^{(\ss,m,k)})$ as follows:
\begin{align*}
\FF^{(\ss,m,k)}&=\{ u\in L^2(I,m+k): u \ll\ss %\text{ is absolutely continuous with respect to } \ss ,\\
 \text{ and } {du\over d\ss} \in
L^2(I,d\ss)\} \\
 \EE^{(\ss,m,k)}(u,v)&=\int_I {du\over d\ss}{dv\over
d\ss} d\ss + \int_I u(x) v(x) k(dx),\ \text{ for }u,v \in \FF^{(\ss,m,k)}.
\end{align*}
It follows from \cite{fhy} that $\FF^{(\ss,m,k)}$ is the closure of the algebra generated by $\ss$
with respect to the norm $\sqrt{\EE^{(\ss,m,k)}(\cdot,\cdot)+(\cdot,\cdot)_m}$. As in \cite{FOT},
if $I=\ang{a_1,a_2}$, we call $a_1$ a regular boundary if $a_1\not\in I$, $\ss(a_1)>-\infty$ and
$m((a_1,c))+k((a_1,c))< \infty$ for some $c\in I$. The regularity of $a_2$ is defined similarly.
Define also
\begin{align*}
\FF_0^{(\ss,m,k)}&= \{u\in \FF^{(\ss,m,k)}: u(a_i)=0 \text{ if } a_i
\text{ is regular boundary   }\};\\
 \EE_0^{(\ss,m,k)}(u,v)&=  \EE^{(\ss,m,k)}(u,v),\ \text{ for }u,v \in
\FF_0^{(\ss,m,k)}.
\end{align*}
When $k=0$, we write it as $(\EE_0^{(\ss,m)},\FF_0^{(\ss,m)})$ for simplicity. The next lemma
asserts that a Dirichlet form is built this way.

\begin{lem} The form
$(\EE_0^{(\ss,m,k)},\FF_0^{(\ss,m,k)})$ is a local irreducible Dirichlet space on $L^2(I;m)$
regular on $I$ and it is strong local if and only if $k=0$.
\end{lem}

\proof We only prove the first statement. The second is clear. Let $J=\ss(I)$ and define a regular
Dirichlet space $(\EE,\FF)$ on $L^2(J,m\circ \ss^{-1})$ (refer to \cite[Example 1.2.2]{FOT} for a
proof) as follows:
\begin{align*}
\FF&=\{ u\in L^2(J,(m+k)\circ \ss^{-1}): u \text{ is absolutely continuous and } u' \in
L^2(J)\} \\
 \EE(u,v)&=\int_J u'(x)v'(x)dx+\int_J u(x)v(x)(k\circ \ss^{-1})(dx),\text{ for }u,v \in \FF.
\end{align*}
Then $(\EE_0^{(\ss,m,k)},\FF_0^{(\ss,m,k)})$ is a state-space transform of $(\EE,\FF)$ induced by
the function $\ss^{-1}$. It shows that $(\EE_0^{(\ss,m,k)},\FF_0^{(\ss,m,k)})$ is a Dirichlet form
on $L^2(I,m)$ by \cite[lemma 3.1]{FFY}. The regularity follows from the fact that $u\circ \ss^{-1}
\in \FF_0^{(\ss,m,k)} \cap C_c(I)$ whenever $u\in \FF\cap C_c(J)$. The local property of
$(\FF_0^{(\ss,m,k)},\EE_0^{(\ss,m,k)})$ is obvious. \qed

\section{Representation of one-dimensional local Dirichlet space}

Let $I$ be an interval or a connected subset of $\R$ and $I^{\circ}$ its interior.

\begin{dfn} \rm A diffusion $X=(X_t,\P^x)$ with life time $\zeta$ on $I$ is a Hunt process on $I$ with
continuous sample paths on $[0,\zeta)$. A diffusion $X$ is called irreducible if for any $x,y\in
I$, $\P^x(T_y<\infty)>0$, where $T_y$ denotes the hitting time of $y$.
\end{dfn}

The irreducibility defined here implies the regularity in \cite{ry} and \cite{rw}. The reason we
use irreducibility is that $I$ is the state space of $X$, while in \cite{ry} and \cite{rw}, $I$ may
contain a trap, thus not a real state space. Another thing which needs to be noted is that a
diffusion defined this way is allowed being `killed' inside $I$, while in some literature it is not
allowed. A diffusion not allowed being killed inside $I$ is called locally conservative. The local
conservativeness is equivalent to the following property: for any $x\in I^{\circ}$, there exist
$a,b\in I$ with $a<b$ and $x\in (a,b)$ such that $\P^x(T_a\wedge T_b<\infty)=1$; if $x$ is the
right (resp. left) end-point of $I$ included in $I$ and finite, then there exists $a\in I$ and
$a<x$ (resp. $a>x$) such that $\P^x(T_a<\infty)=1$. For any regular diffusion $X$, we shall obtain
a process $X'$ through the well-known Ikeda-Nagasawa-Watanabe piecing together procedure. It is
easy to show that $X'$ is a locally conservative regular diffusion on $I$, and $X$ is obtained by
killing $X'$ at a rate given by a PCAF. We say that $X'$ is a resurrected process of $X$ and $X$ is
a subprocess of $X'$.  As VII(3.2) in \cite{ry} or (46.12) in \cite{rw}, a locally conservative
regular diffusion $X$ on $I$ has so-called scale function, namely,
%\begin{thm} Assume that $X$ is a locally conservative irreducible diffusion on $I$. Then
there exists a continuous, strictly increasing function $\ss$ on $I$ such that for any $a,b,x\in I$
with $a<b$ and $a\le x\le b$,
\begin{equation}\label{e:1.1}\P^x(T_b<T_a)={\ss(x)-\ss(a)\over\ss(b)-\ss(a)}.\end{equation} The function $\ss$ is unique
up to a linear transformation.
%\end{thm}
This function $\ss$ is called a scale function of $X$. A diffusion with scale function $\ss(x)=x$
is said to be in natural scale. It is easy to check that if $\ss$ is a scale function of $X$, then
$\ss(X)$ is a diffusion on $\ss(I)$ in natural scale.  A Brownian motion on $I$ is a diffusion on
$I$ which moves like Brownian motion inside $I$ and is reflected at any end-point which is finite
and in $I$ and get absorbed at any end point which is finite but not in $I$. Clearly Brownian
motion on $I$ is clearly in natural scale. Thus Blumenthal-Getoor-Mckean's theorem (Theorem
5.5.1~\cite{BG}) implies that
%\begin{thm}
a diffusion on $I$ in natural scale is identical in law with a time change of Brownian motion on
$I$.
%\end{thm}
More precisely, let $X$ be a locally conservative regular diffusion in natural scale. Then there
exists a measure $\xi$ on $\R$, fully supported on $I$, and a Brownian motion $B=(B_t)$ on $I$ such
that $X$ is equivalent in law to $(B_{\tau_t})$ where $\tau=(\tau_t)$ is the continuous inverse of
the PCAF $A=(A_t)$ of $B$ with Revuz measure $\xi$. The measure $\xi$ is called the speed measure
of $X$. Obviously $X$ is symmetric with respect to $\xi$.

Let now $X$ be an irreducible diffusion on $I$ and $X'$ the resurrected process of $X$ with scale
function $\ss$. Then $\ss(X')$ is symmetric with respect to its speed measure $\xi$ and therefore
$X'$ is symmetric with respect to $\xi\circ\ss$. The diffusion $X$, the subprocess of $X'$, is
certainly still symmetric to $\xi\circ\ss$. An $m$-symmetric Markov process on state space $E$
always determines a Dirichlet form on $L^2(E,m)$. A standard reference for theory of Dirichlet form
is \cite{FOT}, to which we refer for terminologies, notations and results. By results in theory of
Dirichlet form, the Dirichlet form associated with $X'$ is strongly local, irreducible and regular
on $I$. It follows then that the Dirichlet form associated with $X$ is local, irreducible and
regular on $I$. Conversely, given a local, irreducible and regular Dirichlet form on $L^2(I,m)$
with a fully supported Radon measure $m$ on $I$, it is easily seen that the corresponding Markov
process must be an irreducible diffusion on $I$. Therefore one-dimensional irreducible diffusions
are in one-to-one correspondence with one-dimensional local, irreducible and regular Dirichlet
forms. This illustrates that no generality will be lost if we start from such a Dirichlet form as
we shall do in the following sections. In \S2, we shall present a sufficient condition for the
uniqueness of symmetrizing measure. Actually, this condition is almost necessary too. In \S3 we
will give a representation for any 1-dim local, irreducible and regular Dirichlet space. In \S4, we
will give a necessary and sufficient condition for a Dirichlet space to be regular subspace of
another Dirichlet space, which generalizes the main result in \cite{FFY}. As application, two
examples is presented to illustrate that Brownian motion has not only regular extensions and but
also non-conservative regular subspaces.

%We also need to state the following important theorem \cite[Theorem 5.5.1]{BG}. It asserts that if
%two Hunt processes possess identical hitting distribution, then one of them can be transformed to
%the other by a time change.

%\begin{thm}
%Let $X$ and $X'$ be two Hunt processes on $(E,\EE)$ and suppose that, for each compact subset $K$
%of $E_\Delta$, $P_K(x,\cdot)=P'_K(x,\cdot)$ for all $x$. Then there exists a PCAF, $A$ of $X$ which
%is strictly increasing and finite on $[0,\zeta)$ such that if $(\tau_t)$ is the inverse of $A_t$
%and $\hat{X}_t=X_{\tau_t}$ is the time change of $X$, then $\hat{X}$ and $X'$ are equivalent.
% \end{thm}

Fixing an interval $I$ and given a fully-supported Radon measure $m$ on $I$, we shall consider in
this section the representation of a local, irreducible and regular Dirichlet space $(\EE,\FF)$ on
$L^2(I,m)$ in terms of the scale function of the associated diffusion. The form $(\EE,\FF)$ is
assumed to be irreducible, i.e., the associated semigroup is $m$-invariant. Let $X=(X_t,\P^x)$ be
the diffusion process on $I$ associated with $(\EE,\FF)$. It is well known that
%\begin{thm}
the process $X=(X_t,\P^x)$ associated with a local irreducible regular Dirichlet space $(\EE,\FF)$
on $L^2(I,m)$ is an irreducible $m$-symmetric diffusion on $I$. In addition $(\EE,\FF)$ is strong
local if and only if $(X_t,\P^x)$ is locally conservative.
%\end{thm}

Next we give the representation theorem of one-dimensional local, irreducible and regular Dirichlet
space.
\begin{thm}\label{t:1.5}
Let $I=\ang{a_1,a_2}$ be any interval and $m$ a Radon measure on $I$ with $\supp(m)=I$. If
$(\EE,\FF)$ be a local irreducible regular Dirichlet space on $L^2(I,m)$, then
$$(\EE,\FF)=(\EE_0^{(\ss,m,k)},\FF_0^{(\ss,m,k)})$$
where $k$ is a Radon measure on $I$ and $\ss\in \SS(I)$. Furthermore $\ss$ is a scale function for
$(X_t,\P^x)$ which is the diffusion associated with $(\EE,\FF)$.
\end{thm}

%First the following lemma says that two diffusions with the same scale function have the identical
%hitting distributions.
\def\Q{{\mathbf{Q}}}
\proof We shall first assume that $(\EE,\FF)$ is strongly local. Let $\ss$ be a scale function of
$X=(X_t,\P^x)$ associated with $(\EE,\FF)$, and $Y=(Y_t,\Q^x),x\in I$ be the diffusion associated
with Dirichlet space $(\FF_0^{(\ss,m)},\EE_0^{(\ss,m)})$.  Then $X$ and $Y$ have the same scale
function and thus
the same hitting distributions. % By the assumption, we know that $X=(X_t,P^x),x\in I$ and
%$Y=(Y_t,Q^x),x\in I$ possess the same hitting distributions.
It follows from Blumenthal-Getoor-Mckean Theorem that there exists a strictly increasing continuous
additive functional $A_t$ of $X$ such that $(Y_t,\Q^x),x\in I$ and $(\tilde{X}_t,\P^x),x\in I$ are
equivalent, where $\tilde{X}_t=X_{\tau_t}$, and $(\tau_t)$ is the inverse of $(A_t)$.

Note that  $(\tilde{X}_t,\P^x),x\in I$ is $\xi$-symmetric, where $\xi$ is the Revuz measure of $A$
with respect to $m$, and also $m$-symmetric since it is equivalent to $(Y_t,\Q^x),x\in I$. By
Theorem~\ref{t:2.1}, $\xi$ is a multiple of $m$ or $A_t=c t$ for some positive constant $c$. It
shows
that $\tilde{X}_t=X_{t\over c}$. % Let $p_t$, $p'_t$ be respectively the transition functions of
%$(X_t,P^x)$ and $(Y_t,Q^x)$, we have $p_t=p'_{t\over C}$ and
Therefore
\begin{equation*}
\FF=\FF_0^{(\ss,m)} ,\ \EE=c\cdot\EE_0^{(\ss,m)}
\end{equation*}
by (1.3.15) and (1.3.17) in \cite{FOT}.

However scale functions of a linear diffusion could differ by a linear transform. When the scale
function is properly chosen, the constant $c$ above could be 1 (and shall be taken to be 1 in the
sequel). For example $\ss'=\ss/c\in \SS(I)$ is also a scale function for $(X_t,\P^x)$ and we have
\begin{equation*}
\FF=\FF_0^{(\ss',m)} ,\ \EE=\EE_0^{(\ss',m)}.
\end{equation*}

%\noindent\textbf{Proof of Theorem~\ref{t:1.5}:}
In general, when $(\EE,\FF)$ is local, we have the following Beurling-Deny decomposition by
\cite[Theorem 3.2.1]{FOT}
\begin{equation*}
\EE(u,v)=\EE^{(c)}(u,v) + \int_I u(x)v(x) k(dx),\ u,v \in \FF \cap C_0(I),
\end{equation*}
where $\EE^{c}$ is the strongly local part of $\EE$. Define a new symmetric form $(\EE',\FF')$ on
$L^2(I,m+k)$:
 $$\FF'=\FF, \EE'=\EE^{(c)}.$$
 Then $(\EE',\FF')$ is a strongly local irreducible regular Dirichlet space on $L^2(I,m+k)$.
 By the conclusion in the first part, it follows that
 $$\EE^{(c)}=\EE_0^{(\ss,m)},\FF=\FF'=\FF_0^{(\ss,m+k)}=\FF_0^{(\ss,m,k)}.$$ The proof is
completed. \qed

\medskip

\noindent{\bf Remark.} After reading the result above, Professor Fukushima also provides a more
intrinsic proof. Here ``intrinsic'' means a proof without using big theorems developed above but
only using a very profound analysis on the one-dimensional diffusion presented in classical books
K. Ito\cite{I}, \cite{I.57} and Ito-McKean\cite{im}. We shall outline the proof here which is
quoted from Professor Fukushima's e-mail.

\begin{enumerate}

\item Given a diffusion $X$ on a one-dimensional interval $I$, its scale function $\ss$ and speed
measure $m$ are already defined. As you know, $m$ is defined simply by suing the concave property
of the mean exit time from a sub-interval of $I$ when $X$ is locally conservative.

\item ......

\end{enumerate}

%\end{remark}

\section{Regular subspaces}

Let $(\EE',\FF')$ and $(\EE,\FF)$ be two irreducible regular Dirichlet spaces on $L^2(I,m)$. The
space $(\EE',\FF')$ is called a regular subspace of $(\EE,\FF)$ if $\FF'\subset\FF$ and
$\EE(u,v)=\EE'(u,v)$ for any $u,v\in\FF'$. All non-trivial regular subspaces of linear Brownian
motion is characterized clearly in \cite{FFY}. In this section we shall further give a necessary
and sufficient condition for $(\EE',\FF')$ to be a regular Dirichlet subspace of $(\EE,\FF)$,
which extends the result in \cite{FFY}. %A
%natural problem is whether $(\EE',\FF')$ is a regular Dirichlet subspace of $(\EE,\FF)$. In this
%section we shall give a necessary and sufficient condition .

%Denote by $\AA$ the totality of measurable subsets $A$ of $I$ satisfying that, for any $x,y\in I$
%with $x<y$, $|(I\setminus A)\cap (x,y)|>0$, i.e., the complement of $A$ has a positive measure on
%any non-empty open subinterval. Two sets in $\AA$ are regarded to be equivalent if they differ by a
%zero-measure set.
Using the representation in \S3, we have
\begin{align*}
(\EE,\FF)&=(\EE_0^{(\ss_1,m,k_1)},\FF_0^{(\ss_1,m,k_1)});\\
(\EE',\FF')&=(\EE_0^{(\ss_2,m,k_2)},\FF_0^{(\ss_2,m,k_2)}),
\end{align*}
where $\ss_1,\ss_2\in \SS(I)$ and $k_1$, $k_2$ are two Radon measures on $I$. Now comes our main
result.

\begin{thm}
Let $(\EE',\FF')$ and $(\EE,\FF)$ be two local irreducible regular Dirichlet spaces on $L^2(I,m)$.
Then $(\EE',\FF')$ is a regular subspace of $(\EE,\FF)$ if and only if
 \begin{itemize}
 \item[$(1)$] $k_1=k_2$,
 \item[$(2)$] $d\ss_2$ is
absolutely continuous with respect to $d\ss_1$ and the density %there exists a set $A\in \AA$ such that
${d\ss_2 / d\ss_1}$% =  1_{A^c}.$$
is either $1$ or $0$ a.e. $d\ss_1$.
 \end{itemize}
\end{thm}
\proof \
  It suffices to prove it for the case that both $(\EE',\FF')$ and $(\EE,\FF)$
are strongly local. Assume that $\FF' \subseteq \FF$ and let $(X_t,\P_x)$ and $(X'_t,\P'_x)$ be the
diffusion processes associated with $(\EE,\FF)$ and $(\EE',\FF')$, respectively. For any
$a<c<x_0<d<b$, define $$u_{\{c,d\}}^{x_0}(x):=\P'_x(T_{x_0}<T_{\{c,d\}}).$$ We have
$u_{\{c,d\}}^{x_0}(x)\in \FF' \subseteq \FF$, and  it shows that $u_{\{c,d\}}^{x_0}(x)$ is
absolutely continuous with respect to $\ss_1$, while $u_{\{c,d\}}^{x_0}$ is a linear transformation
of $\ss_2$ on $(c,x_0)$. It follows that $d\ss_2$ is absolutely continuous with respect to $d\ss_1$
on $(c,x_0)$. Similarly it is also true on $(x_0,d)$. Taking $(c,d) \uparrow (a,b)$, it follows
that $d\ss_2$ is absolutely continuous with respect to $d\ss_1$. Let $f:={d\ss_2 / d\ss_1}$. Then
we have
\begin{align*}
\EE'(u,v)&=\int_I {du\over d\ss_2}{dv\over d\ss_2}d\ss_2; \\
\EE(u,v)&=\int_I {du\over d\ss_1}{dv\over d\ss_1}d\ss_1 \\
       &=\int_I {du\over d\ss_2}{dv\over d\ss_2}f^2 d\ss_1\\
       &=\int_I {du\over d\ss_2}{dv\over d\ss_2}f d\ss_2
\end{align*}
for any $u,v \in \FF'$. It follows then that $fd\ss_1=f^2d\ss_1$ and that either $f=0$ or $f=1$
a.e. with respect to $d\ss_1$. Since $\ss_1$ and $\ss_2$ are continuous and strictly increasing,
$f$ has the property that for any $x,y\in I$ with $x<y$, \begin{equation}\int_x^y
1_{\{f=1\}}d\ss_1>0.\end{equation}
 The converse is obvious from the above discussion.
\qed

%The following theorem is an easy consequence.
Let now $$(\EE,\FF)=(\EE_0^{(\ss,m,k)},\FF_0^{(\ss,m,k)})$$ be a local irreducible regular
Dirichlet spaces on $L^2(I,m)$. Take a Borel set $A$ having property that for any $x,y\in I$ with
$x<y$, \begin{equation}\label{e:3.2}\int_x^y 1_{A^c}d\ss>0.\end{equation} Define
$d\ss_0=1_{A^c}\cdot d\ss$. Then $\ss_0\in \SS(I)$ and $(\EE_0^{(\ss_0,m,k)},\FF_0^{(\ss_0,m,k)})$
is a regular subspace of $(\EE,\FF)$. It is easy to check that
$$\FF_0^{(\ss_0,m,k)}=\{u\in\FF: \ du/d\ss=0\ \text{a.e. with respect to}\ d\ss\ \text{on}\ A\}.$$
Hence we have a corollary.

\begin{C}
For any Borel set $A$ satisfying \text{\rm (\ref{e:3.2})},
\begin{equation}
\FF^{A}=\{u\in\FF: \ du/d\ss=0\ \text{a.e. with respect to}\ d\ss\ \text{on}\ A\}
\end{equation}
is a regular subspace of $(\EE,\FF)$. Conversely any regular subspace of $(\EE,\FF)$ is induced by
such a set.
\end{C}

Finally, we shall give two interesting examples. The first example is a local irreducible and
regular Dirichlet space which takes the Dirichlet space $(H^1([0,1]),{1\over 2}\DD)$ of reflected
Brownian motion on $[0,1]$ as a proper regular subspace.
\\
%We also can observe that for this series, neither {\it largest} nor {\it smallest} Dirichlet space exists. \\

\noindent {\bf Example 1.} Let $c(x)$ be the standard Cantor function on $[0,1]$ and let
$\ss(x):=x+c(x)$. Take $m$ to be the Lebesgue measure on $[0,1]$. Then the Dirichlet space
$(H^1([0,1]),{1\over 2}\DD)$, corresponding to Brownian motion on $[0,1]$, is a regular subspace of
$(\FF^{(\ss,m)},{1\over 2} \EE^{(\ss,m)})$ by the theorem above and $H^1([0,1])$ is properly
contained in $\FF^{(\ss,m)}$.\\

The second example shows that 1-dim Brownian motion has a non-conservative regular subspace. For
this we state a criterion for irreducible one-dimensional diffusions to be conservative (see
\cite{rw}). %The result is essentially classical and we include a proof here since we use the
%approach of Dirichlet forms.
%Let $I$ be a linear connected set and $m$ be a Radon measure on $I$.
Let
$$(\EE,\FF)=(\EE_0^{(\ss,m,k)},\FF_0^{(\ss,m,k)})$$
where $k$ is a Radon measure on $I$ and $\ss\in \SS(I)$, be a local, irreducible and regular
Dirichlet space on $L^2(I,m)$ and $X=(X_t,P_x)$ the associated diffusion. In this case it is either
recurrent or transient. We call the left endpoint $a$ of $I$ is
\begin{itemize}
\item[\rm (1)] of the first class if $a$ is finite and $a\in I$;
\item[\rm (2)] of the second class if $a\not\in I$ and
$\ss(a)=-\infty$; \item[\rm (3)] of the third class if $a\not\in I$ and $\ss(a)>-\infty$.
\end{itemize}
%Suppose that $a$ is the left endpoint of $I$ and $\ss$ is a scale function.
We call $a$ is dissipative if $a$ is of the third class and
\begin{equation}\label{4.1}
\int_a^c (\ss(x)-\ss(a)) m(dx)<\infty
\end{equation}
for some $c\in I$, and hence for all $c\in I$. Obviously, the finiteness (\ref{4.1}) is independent
of the choice of the scale function $\ss$ and the point $c$. If $a$ is not dissipative, we call it
conservative. The dissipativeness and conservativeness for the right endpoint may be defined
similarly. Fix a point $c>a$, define $M(x):=m((x,c))$ for $a<x<c$.

\begin{lem}
The left end-point $a$ is dissipative if and only if $a$ is of the third class and
\begin{equation}\label{4.2}
\int_a^c M(x) d\ss(x)<\infty.
\end{equation}
If $a$ is dissipative, $\lim_{x\downarrow a} M(x)\ss(x)=0$. Similar conclusions hold for the right
end-point.
\end{lem}

%\proof \ Firstly, we prove the second assertion. Assumed that $\overline{\lim}_{x\downarrow a}
%M(x)\ss(x)=p>0$. Then there exists a sequence $\{r_n\} \leq c$ such that $r_n \downarrow a$ and
%$M(r_n)\ss(r_n)>{p\over 2}$, $\ss(r_{n+1})<{1\over 2}\ss(r_n)$. Therefore,
%\begin{eqnarray*}
%\int_a^c M(x)d\ss(x)&\geq&\sum_{n}\int_{r_{n+1}}^{r_n} M(x)d\ss(x) \\
%&\geq& \sum_{n} M(r_n)(\ss(r_n)-\ss(r_{n+1})) \\
%&\geq& \sum_{n} {p\over 4} =\infty.
%\end{eqnarray*}
%It contradicts with (\ref{4.2}). Similarly, we also can prove that (\ref{4.1}) implies that
%\newline $\lim_{x\downarrow a} M(x)\ss(x)=0$.

%Using integration by parts,
%\begin{eqnarray*}
%\int_a^c M(x)d\ss(x)&=&\lim_{a'\downarrow a} \int_{a'}^c M(x)d\ss(x)\\
%&=& \lim_{a'\downarrow a} (M(c)\ss(c)-M(a')\ss(a')-\int_{a'}^c
%\ss(x) d M(x))\\
%&=& M(c)\ss(c)-\lim_{a'\downarrow a}M(a')\ss(a')-\lim_{a'\downarrow a} \int_{a'}^c \ss(x)m(dx)\\
%&=& M(c)\ss(c)-\int_a^c \ss(x)m(dx).
%\end{eqnarray*}
%It completes the proof. \qed\\

%Next, we shall prove sufficient and necessary conditions for one-dimensional diffusion to be
%transient, recurrent or conservative.

%\begin{thm}

\begin{thm} The Dirichlet space $(\EE,\FF)$ (or $X$) is
\begin{itemize}
\item[\rm (1)] recurrent if and only if $k=0$ and both endpoints
are of the first class or the second class;
%\item[\rm (2)] transient if and only if it is not recurrent;
\item[\rm (2)] conservative if and only if $k=0$ and both endpoints are conservative.
\end{itemize}
\end{thm}

We now give an example which illustrates that the Dirichlet space $({1\over 2}\DD,H_0^1(\R))$ of
Brownian motion on the real line $\R$ has non-conservative regular subspaces, comparing an example
in \cite{FFY} which shows Brownian motion has transient regular subspaces.
\\

 \noindent {\bf Example 2.} Define a local irreducible and regular
 Dirichlet space $(\EE_0^{(\ss,m)},\FF_0^{(\ss,m)})$ on $L^2(\R,m)$, where $m$ is the usual Lebesgue
 measure, by giving a scale function
\begin{equation*}
\ss(x) =\int_0^x1_G(y)dy,\ x\in \R,
\end{equation*}
where
\begin{equation}\label{4.3}
G = \bigcup_{r_n\in Q}\left(r_n-\frac{1}{2^{n+1}},r_n+\frac{1}{2^{n+1}}\right),
\end{equation}
where $Q$ is the set of positive rational numbers. We choose an order on $Q$ as follows: if $a,b\in
Q$, and $a={q_1\over p_1}$, $b={q_2\over p_2}$ take the simplest form, we define
$$a \prec b \Leftrightarrow {\rm either}\ p_1+q_1 <p_2+q_2 {\rm \ or\ } p_1+q_1 =p_2+q_2{\rm \ and\ }q_1<q_2
.$$ Then the order $\prec$ makes $Q$ a sequence $\{r_n\}$ in (\ref{4.3}). Clearly $r_n\leq n$. Thus
\begin{align*}
\int_0^\infty x d\ss(x) \le\sum_n \int_{(r_n-\frac{1}{2^{n+1}},r_n+\frac{1}{2^{n+1}})} x dx =\sum_n
{r_n\over 2^n}  \leq \sum_n {n\over 2^n} <\infty.
\end{align*}
This shows the right endpoint is dissipative. Therefore the associated process is not conservative.
\\

\noindent \bf Acknowledgements: \rm The authors would like to thank Professor M.Fukushima for his
helpful suggestions.

%\vskip 0.6truein

\noindent Addresses:

\medskip
\noindent {\bf X. Fang}: Department of Mathematics, Fudan University, Shanghai, China.

Email: {\texttt fangxing@fudan.edu.cn}

\medskip

\noindent {\bf J. Ying}: Department of Mathematics, Fudan University, Shanghai, China.

Email: {\texttt jgying@fudan.edu.cn}

\medskip

\noindent {\bf M. Zhao}: Department of Mathematics, Zhejiang University, Hangzhou, China

e-mail: {\texttt zhaomz@zju.edu.cn}

\end{document}